# Diverse yet consistent: How mathematicians position computational thinking across research and teaching

Jan-Fredrik Olsen[1,2] (jan-fredrik.olsen@math.lu.se),
Tor Ole B Odden[2] and Elise Lockwood[3]

**Abstract:** Recent research in mathematics education points to an “epistemic clash” when programming and computational thinking (CT) are leveraged alongside more established forms of mathematical thinking (MT). The emergence of generative AI emphasises the need to understand the mechanisms shaping relations between CT and MT.

We address this need by analysing interviews with 15 mathematicians on their use of computations across their teaching and research activities. The interviews were conducted at a critical site with a history of integrating computations across its science and mathematics programs for more than 20 years. Drawing on Cultural Historical Activity Theory and Communities of Practice theory, we consider MT and CT as methodologies grounded in practice.

We identify three perspectives shaping how mathematicians position CT: mathematical theory considered as a source of control, computations as a source of pragmatic reach, and real-world impact as a source of legitimacy. This three-perspectives model explains why mathematicians who emphasise real-world impact are most likely to carry programming into teaching, whereas those who position theoretical mathematics as authoritative are least likely to do so. Mathematicians working on numerical algorithms occupy an uneasy intermediate position.

Our findings suggest that the perceived clash between MT and CT is not purely epistemic, but also ontological, as it depends on how computations are positioned within the goal of doing mathematics. For mathematics education, this implies that perceived meaningful integration with CT is mediated by context, and that more extensive use can be stabilised by leveraging authentic learning goals external to mathematics.

[1] Lund University, Lund, Sweden.
[2] University of Oslo, Oslo, Norway.
[3] Oregon State University, Corvallis, OR, USA

# Introduction

Throughout its history, mathematics has been caught in a balancing act between an internal identity as a timeless language for understanding nature and an external appearance as a tool for bookkeeping and computation. This balancing act between pragmatic use-value and theoretical authority persists across mathematical practices today, and is at the same time a cause of tension and a driver of innovation. While tools such as Wolfram Alpha, MATLAB, Maple, Python, and generative AI are regularly used by mathematicians and students to check examples and conjectures, their legitimacy in mathematics education is contested (e.g., Mørken & Lockwood, 2019).

Recently, educational researchers point to signs of an epistemic clash between the use of programming and mathematical 'registers' in the classroom (Bråting & Kilhamn, 2021; Borg & Fahlgren, 2025; Huang et al., 2025). They suggest that students find it hard to coordinate conceptual mathematical work with computational work, which shifts their focus to debugging code implemented line by line. These findings contrast with the work of Papert (1980) and diSessa (2000), who advocate a view of computation as affording a more playful approach to learning mathematics. The rise of generative AI risks exacerbating such tensions. Where some see a tool that allows us to rethink and strengthen the interface between mathematical and computational work in both research and teaching, others see a threat to scientific and educational engagement and integrity (Nguyen & Pham, 2025).

From a historical perspective, mathematical research, rather than teaching, appears to be the main arena for innovation and the adoption of new ideas. In light of the current rapid pace of technological development, the following motivating questions arise: how are mathematicians integrating computational technologies into their research, and how does this affect how they position computations as a tool for learning in their own teaching?

In this paper, we address these questions by analysing interviews with 15 mathematicians situated at a university where the integration of computations into education has been institutionally supported over the last 20 years. We frame our study in an interpretive tradition and draw on Cultural Historical Activity Theory (Engeström, 2015) and Communities of Practice theory (Lave & Wenger, 1991; Wenger, 1998) in our analysis of how mathematicians position computations across their research and teaching practices.

We argue that mathematicians' positioning of computations in teaching is shaped by how they configure three perspectives within their research: mathematical theory as a source of control, computations as a source of pragmatic reach, and real-world relevance as a source of legitimacy. Extending our previous work (Olsen et al., 2016), we develop a three-perspectives model that explains how different configurations of the research object pattern mathematicians' engagement with computations across research and teaching. This allows us to show that perceived tensions and synergies between mathematical and computational thinking are not only epistemic, but also depend on how computations are positioned in relation to what counts as doing mathematics.

# Background

## Computational thinking in educational research

Although computational thinking (CT) has a long history in mathematics education research, there is no consensus on a definition (see, e.g., Gadanidis et al., 2021). We adopt a stance similar to Brennan and Resnick (2012) and Weintrop et al. (2016), who conceptualise CT as methodologies grounded in professional practice, and not just as thinking skills. We take an analogous view of mathematical thinking (MT), understanding it as grounded in the professional practices of mathematicians. In a consensus-driven Delphi study, Kallia et al. (2021) found that mathematicians and mathematics teachers view CT as a "*structured problem-solving approach"* that leverages techniques such as algorithmic thinking, abstraction, and pattern recognition. This suggests that mathematicians' conceptualisation of CT lies closer to MT than to Wing's (2006) well-known framing of CT as a set of transferable skills from computer science.

Huang et al. (2021) found that mathematics teachers tend to adapt their articulate conception of CT to the context. In mathematics, they described CT in ways that aligned with familiar forms of mathematical problem solving and reasoning, whereas in computing, their descriptions were strongly influenced by programming. Similarly, Bråting and Kilhamn (2021) argue that programming in mathematics can foreground process and time-related change while backgrounding relational meanings central to algebra (see also Huang et al., 2025; Borg & Fahlgren, 2025). Taken together, this line of work resonates with what Huang et al. (2025) describe as a "clash of epistemic tools", suggesting that how tasks are framed and which tools are brought into play can shift how CT and MT are positioned relative to one another within mathematical activity.

Taken together, this research suggests that tensions between the use of CT and MT may depend less on abstract definitions of either than on how they are positioned relative to one another within mathematical activity. This, in turn, suggests that these positionings can be better understood by studying how they shape mathematicians' professional mathematical practice.

## Computational thinking in mathematicians' practices

One way to gain insight into how CT and MT can be positioned in mathematical work is to examine mathematicians' professional practices. However, while such practices have long been considered a valuable focus of mathematics education research, little research has been done on mathematicians' use of computations (Mørken & Lockwood, 2019). Here, we review some of the research done on mathematicians' practices and how it relates to the use of computations and computational thinking.

Within mathematics education research, a consistent picture has emerged of mathematical research practice being more heterogeneous than their teaching practices alone might suggest (Burton, 2004). This has led some educational researchers to warn about methodological difficulties related to interpreting such work and its implications for education (Weber et al., 2020).

This picture of heterogeneous practices carries over to findings on how mathematicians use computations. In two survey studies, Lavicza (2010) and Buteau et al. (2014) reported high use of technology among mathematicians across research and teaching. However, this use was uneven across practices; in the latter study, mathematicians reported significantly higher programming use in their research than in their teaching (43% versus 18%). In the UK context, Sangwin and O'Toole (2017) found that programming tasks were invisible on assessment in mathematics courses, despite programming being a compulsory part of their modules.

We are aware of only a small number of interview studies on mathematicians' use of computations in their professional practices, and these offer limited insight into the epistemic reasons why such use appears to be less prominent in teaching. Lockwood et al. (2019) argue that computing is part of mathematicians' disciplinary practice, highlighting uses such as experimentation, conjecturing, and the construction of their own tools, but they do not focus on how such practices carry over into teaching. Broley et al. (2018), observing similar practices, propose a taxonomy based on levels of engagement with programming to describe variation in how mathematicians use code. They report that mathematicians point to institutional constraints (curricular, departmental and cultural) to explain why the use of computations does not carry over to teaching.

In our recent work (Olsen et al., 2026), we observe that what at first appeared to be heterogeneous accounts of the role of computations in mathematics could be organised according to a small number of coherent *epistemic archetypes* based on their preferred research orientation (pure or applied) *and* use of computations. This resulted in four archetypes, each with its internally coherent stories of how the legitimacy of computations was judged across teaching and research. However, this did not mean their actual engagement was the same; rather, it meant their ways of reasoning were shared.

In light of the above research, the question of an epistemic clash between CT and MT in mathematicians' practices seems to be primarily about how these methodologies are conceptualised and positioned relative to each other, and less about whether CT can support MT. Addressing this issue is particularly important as a clearer understanding of such interfaces may provide a basis for further efforts to examine how computations, programming, computational thinking, and emerging technologies can be productively leveraged in the context of mathematics education. By addressing the fundamental mechanisms at play, we hope to also shed light on how we can understand and shape the emerging role of generative AI in mathematics education.

Specifically, we seek to investigate how mathematicians' research orientation shapes their positioning of computations and CT as tools. Our research questions are:

- **RQ1:** How do mathematicians position computational thinking relative to mathematical thinking within their research activity?
- **RQ2:** How does this positioning mediate engagement with computational thinking in teaching?
- **RQ3:** How do these positionings relate to mathematicians' talk about tensions or synergies between mathematical and computational thinking *in education*?

# Theoretical framework

We combine Engeström's Cultural Historical Activity Theory (CHAT) (Engeström, 2015) and Lave and Wenger's Community of Practice theory (CoP) (Lave & Wenger, 1991; Wenger, 1998) to analyse how mathematicians position computation within their professional practice. CHAT foregrounds the object structure and tensions of research and teaching activity, while CoP helps us interpret boundaries within and across specialised mathematical communities.

## Cultural Historical Activity Theory (CHAT)

Drawing on CHAT (Engeström, 2015), we conceptualise mathematicians' research and teaching practices as a situated, collective, enduring, and goal-oriented activity. We treat activity as an analytic construct that contains overlapping sub-activities. In our context, the overall activity is what the interviewed mathematicians do in a professional capacity, which contains the more specialised activities of teaching and research (e.g., Núñez, 2009; Jaworski, 2010; Fredriksen & Hadjerrouit, 2020).

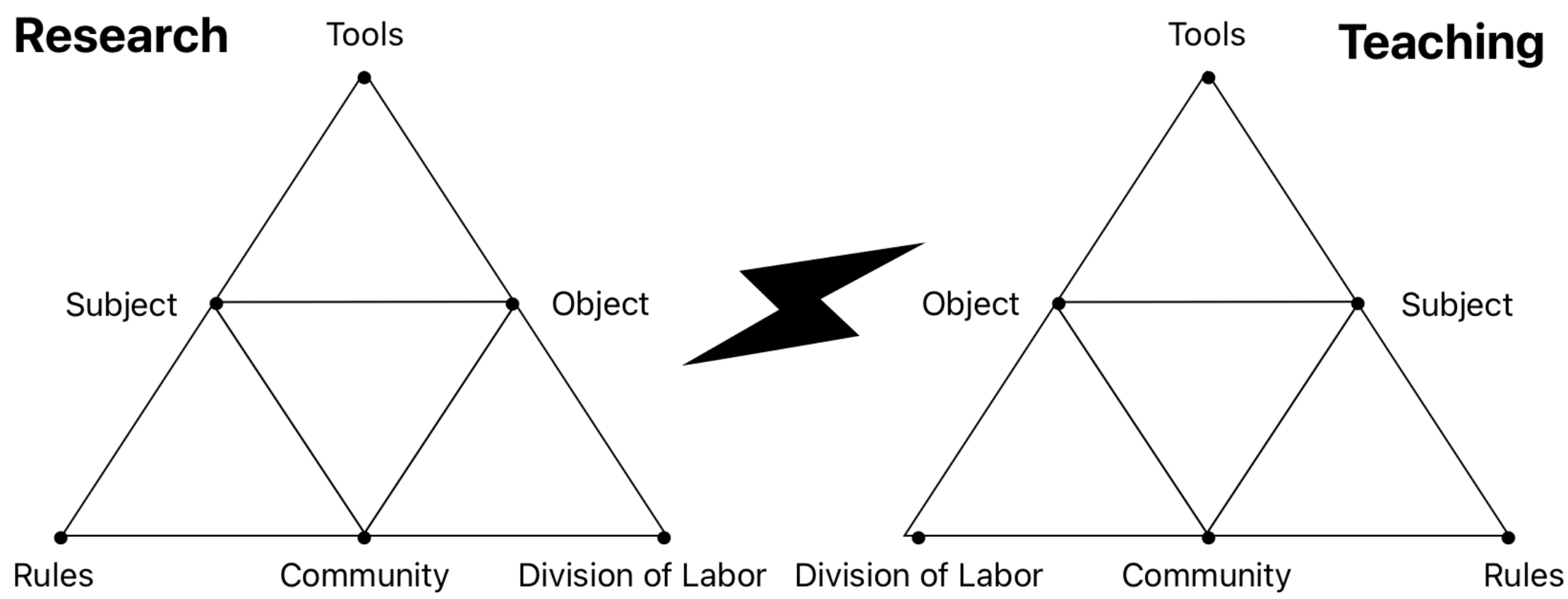

Figure 1: The unit of analysis of Engeström's CHAT is usually taken to be a pair of activities that interact. In our case, we can consider the activities of "mathematics research" and "mathematics teaching" as sub-activities of "working as a mathematician".

Engeström conceptualises activity as an extension of Vygotsky's mediating triangle (Figure 1). The *object* is the collective purpose of the activity, serves as its centre of gravity and main organising principle (e.g., new mathematical insight), and can differ from the concrete, physical *outcome* of the activity (e.g., theorems or research publications). *Subjects* (mathematicians), situated in a *community,* work toward this object by leveraging *tools* (material and immaterial). Subjects are driven by personal *motives* that may or may not be well-aligned with the object. Their actions are shaped by *rules* (formal and informal) and a *division of labour* (formal and informal).

Our central concern is the relation between object and tool. More specifically, we examine how mathematicians position computations and CT in relation to the purpose of their mathematical work and MT: whether they are seen as constitutive of that purpose, or mainly as means for achieving it. Because the object is an analytic construct, its identification is necessarily interpretive and depends on the scope at which the activity is framed. In this study, the object is inferred by triangulating participants' accounts of what they see themselves contributing to collectively through research and teaching, and how they position computations in relation to that contribution.

When serving as an analytic framework, CHAT directs our attention to tensions and contradictions that become visible when studying an activity from vantage points provided by the nodes of the activity triangle (Engeström & Mwanza, 2005; Olsen et al., 2016). These tensions are considered manifestations of constitutive, historical and irremovable forces that shape the development of the activity (Engeström, 2015; Roth & Radford, 2011). In our context, a central dialectical contradiction concerns the boundary between what counts as mathematics and what does not. This contradiction may manifest within or between nodes of the activity system, for example, in relation to the object (what the activity is directed toward), the rules (what is permitted), and the tools (how the activity may be carried out).

## Communities of practice

To complement the systemic perspective of CHAT, we draw on Lave and Wenger's Community of Practice theory (CoP) (Lave & Wenger, 1991; Wenger, 1998). In CoP, emphasis is placed on the lived experience of subjects in a community of professional practice, and concepts serve to link their accounts of that experience to a professional landscape formed by mutual engagement, shared repertoire, and joint enterprise.

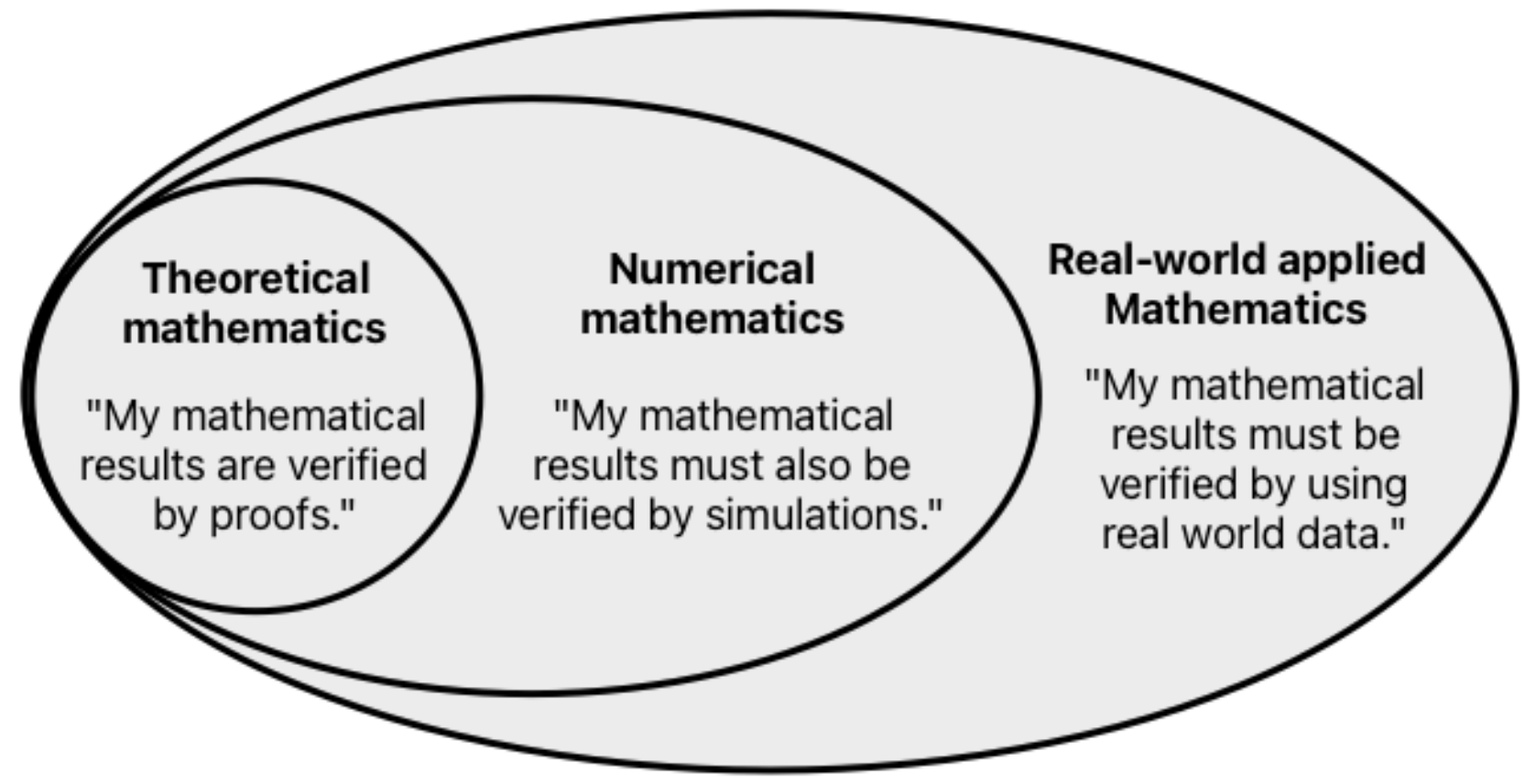

Figure 2: We can conceptualise mathematics as consisting of communities that interact along boundaries of practice. When we cross borders, meanings may change, leading to conflict and innovation.

In our study, CoP provides us with concepts that allow us to discuss boundaries within and between communities of mathematicians. Although a mathematics department can be framed as a community of practice, its sense of mutual engagement, shared repertoire, and joint enterprise may seem stronger from outside than from the inside. Indeed, mathematical work is highly specialised, with engagement taking place mostly within smaller research groups. Such close and sustained interaction leads to specialisations, with diverging repertoires, terminologies and understandings of what should count as mathematics (Figure 2). This creates boundaries and barriers within the community, leading to questions of alignment and positioning; which practices count as central or peripheral, legitimate or marginalised?

Combining the perspectives of CHAT and CoP, we obtain a framework that allows us to consider the activities of mathematicians situated in a community characterised by boundaries arising from high levels of specialisation. Activity close to the boundary of a community is associated with friction and conflict. However, it also requires negotiating meaning, which can lead to innovation and new ideas. Changes in the object, rules, tools, community, and division of labour characterise the boundaries between activities.

# Methodology

We conducted a qualitative interview study of mathematicians at a critical site using CHAT and CoP to analyse how they position CT across research and teaching (Jaworski, 2010; Olsen et al., 2026).

## Context

We conducted the study at a major Northern European university, where the integration of numerical computations into mathematics and science study programs has been an officially sanctioned policy for over 20 years. We treated the site as critical because long-standing institutional support for computation reduced infrastructural and curricular constraints, allowing epistemic issues to come more clearly into view (Flyvbjerg, 2006).

Because mathematicians with low programming engagement might not see the study as relevant, we used sequential purposive sampling with referral to maximise variation in research specialisation and programming engagement. The final sample comprised 15 mathematicians from all formal research units in the department and included early-, mid-, and late-career researchers. Most participants were men, reflecting gender imbalance at the site and a greater willingness of men to participate.

## Data collection

We conducted 15 semi-structured interviews between September and October 2023, each lasting 1.5-2 hours. Participants were first asked about research practices and uses of computation, then about teaching, and finally about institutional context. Interviews were audio-recorded, transcribed, and anonymised. Because only one interview was conducted in English, most excerpts were translated. The interviews were conducted by the first author, whose background in theoretical mathematics enabled participants to use disciplinary language freely and supported detailed probing of mathematical content.

## Operationalisation of the theoretical framework

Following Mwanza and Engeström (2005), we operationalised key CHAT and CoP constructs as questions for analysing how participants positioned CT in relation to their activity (Table 1).

In this study, we conceptualise MT and CT as situated, practice-based categories. These categories are therefore inferred from participants' accounts of professional practice rather than treated as cognitive thinking skills.

| Table 1: Operationalisations of theoretical constructs. | |
|---|---|
| **Theoretical construct** | **Operationalisation** |
| Object | What does the mathematician consider to be the overall purpose or goal of doing mathematical research within their community?<br><br>*“And that is what I would call insight. Yes, that is insight - and fundamentally, that is what I am after.*” |
| Tool | What material and immaterial tools does the mathematician consider to be of help in their research activity?<br><br>“*I only do simple things - I only use simple numerical methods - so it is not any advanced program with several thousand lines…*” |
| Central positioning | What does the mathematician consider to be a natural and valuable part of their activity?<br><br>*“No, if I had not had computers and programming as tools, I would not have managed to accomplish anything.*” |
| Peripheral/marginal positioning | What does the mathematician consider to be auxiliary parts of the activity that are “nice to haves” but that could be dispensed with?<br><br>“*... and implementation of these methods is not the most important part of my research, you could say.*” |
| Internal and external positioning | What does the mathematician see as part of their activity, and what do they not see as part of their activity?<br><br>“*I do not see the use of programming to implement an algorithm as part of mathematics.*” |
| Mathematical thinking (MT) | What do mathematicians perceive as tools (material and immaterial) for solving problems by theoretical means? |
| Computational thinking (CT) | What do mathematicians perceive as tools (material and immaterial) for solving problems by computational means? |
| Real-world thinking (RW) | A catch-all term we use for the shared repertoires of communities in the empirical sciences, such as physics, chemistry, and biology. |

## Analysis

Analysis was iterative and based on coding, case summaries, and memoing informed by CHAT and CoP. To reduce risks of overalignment and tacit misunderstanding, we treated interviews as the unit of analysis and interview meaning as co-constructed by the interviewer and participant rather than as transparent reports of practice. The first author first conducted a framework analysis of a subset of interviews and then developed case summaries for all participants' research and teaching practices. Comparison across cases supported the development of themes around programming in research and teaching and, ultimately, the three-perspectives model used in the final analysis (Figure 3, Table 2).

Classifications were based primarily on interview accounts, supported by publication lists and institutional homepages. Most participants aligned clearly with the model; ambiguous cases were resolved through continued comparison as the model was refined. Where participants described trajectories of change, classification reflected their current position. Interpretations were discussed within the interdisciplinary author team and later presented to a workshop of mathematicians to assess plausibility and explanatory value. To preserve anonymity, we deemphasise individual portraits and do not use pseudonyms.

Generative AI was used by the first author to polish language and condense certain paragraphs, with significant human oversight.

# Results

## A three-perspectives model of CT in the research object

Our analysis identified three perspectives that patterned how mathematicians positioned CT in relation to their research object:

- **The mathematical thinking perspective (MT)** concerns the extent to which theory is positioned as a source of epistemic control.
- **The computational thinking perspective (CT)** concerns the extent to which computations are positioned as external to, supportive of, or constitutive of the research object.
- **The real-world thinking perspective (RW)** concerns the extent to which applications function mainly as motivation or also as a source of engagement and legitimacy.

Using these perspectives as an analytic lens, we identified three predominant configurations among the interviewed mathematicians: MT, MT-CT, and MT-CT-RW (Table 2, Figure 3).

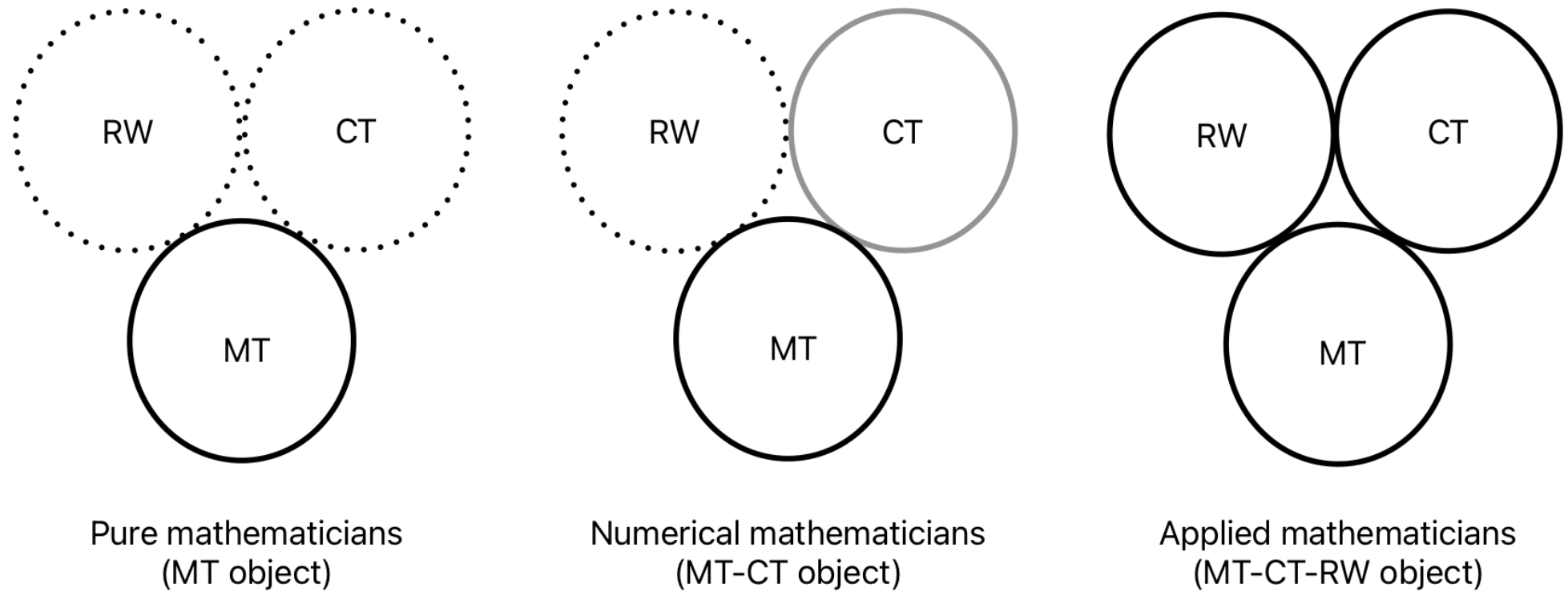

Figure 3: Observed orientations of the research objects represented by the interviewed mathematicians. Solid lines indicate overt and legitimate recognition of methodology as an internal part of their research object, while dotted lines indicate external positioning. Grey lines indicate internal but subordinate positioning.

These groups are analytic categories grounded in participants' orientations across the three perspectives, and the labels "pure mathematicians," "numerical mathematicians," and "applied mathematicians" are intended as descriptive rather than normative. Boundary cases were classified by the source of epistemic legitimacy that participants currently foreground.

Table 2: Summary of positionings of CT in research and teaching by composition of MT, CT, and RW in the object. Group labels are descriptive and grounded in our situated data.

| **Object configuration** | **Positioning of CT in research** | **Positioning of CT in teaching** |
|---|---|---|
| **Pure mathematicians**<br><br>Object: MT<br><br>(n=6) | CT external to the object of research. Positioned as supportive but non-constitutive.<br><br>CT is positioned external to the object. As a tool, CT is epistemically supportive but non-constitutive.<br><br>RW is positioned externally as a source of motivation, including when seen as mediating the relevance of the interface between CT and MT. | Use of CT carries across from research when seen as low-threshold motivational, vocational, and epistemic support.<br><br>Emphasis on RW can increase use of CT in teaching, but such use tends to be viewed as external and vocationally motivated, and seen as a distraction.<br><br>Use of CT increases at the thesis level. |
| **Numerical mathematicians**<br><br>Object: MT-CT<br><br>(n=3) | The object of research is the discovery and validation of theoretical insight into efficient algorithms for solving mathematical models.<br><br>CT is positioned as internal and constitutive to the object. However, as a tool, CT is subordinate to MT and avoided when the algorithmic nature of theoretical concepts has low visibility.<br><br>RW is positioned externally as a source of motivation, including when seen as mediating the relevance of the interface between CT and MT. | CT is more likely to carry over when the object of teaching is considered to have a computational nature. For concepts seen as purely theoretical, the use of CT is suppressed and viewed as a distraction.<br><br>Emphasis on RW can increase use of CT in teaching, but such use tends to be viewed as external and vocationally motivated, and seen as a distraction. |
| **Real-world applied mathematicians**<br><br>Object: MT-CT-RW<br><br>(n=6) | The object of research is the discovery, construction, validation, and deployment of efficient algorithms for solving mathematical models.<br><br>CT, RW, and MT are all considered internal and constitutive to the object. However, the emphasis is uneven, and we observe three main variants (method, model, or implementation emphasis).<br><br>All variations view MT as a source of epistemic control and refer to RW as a source of legitimacy for the use of CT. | While RW stabilises the relevance of CT in teaching, CT carries over differently across the three different emphases.<br><br>Method emphasis: Programming use carries over most clearly, with most seeing 'simulations' as part of pedagogy. Emphasis on MT as a source of control can temper this.<br><br>Model emphasis: The strong emphasis on MT as a source of control can suppress the use of CT in foundational courses. Use converges to that in research at higher levels.<br><br>Implementation emphasis: A pragmatic focus on MT-CT as tools can lead to CT being used only when necessary. This leads to suppression of CT in foundational courses. |

## Pure mathematicians: CT as external support

We consider six of the interviewed mathematicians to be oriented toward classical theoretical work. We interpret their object of research as the discovery of general relations between abstract concepts. In their research, we see a tension between positioning CT as an external distraction versus as a mediator of MT insight (see Figure 4).

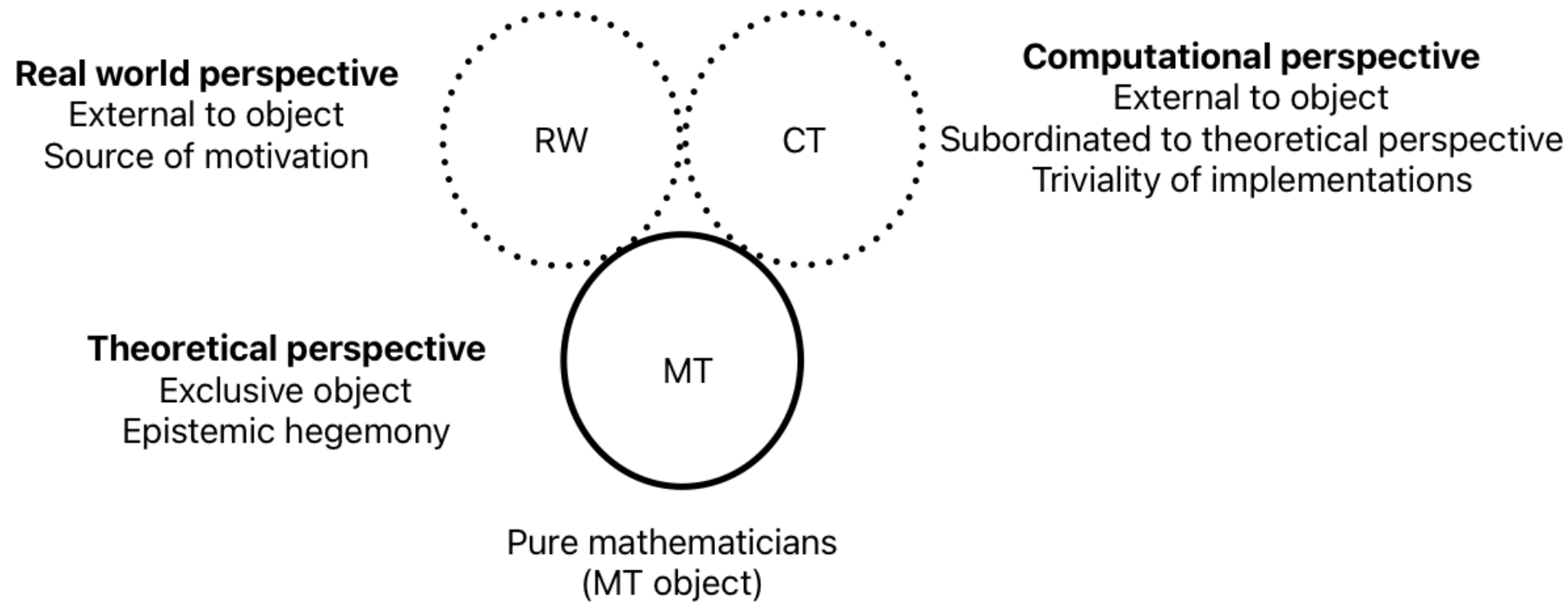

Figure 4: The purely theoretical mathematicians conceptualise their object as being only about MT, and position CT as external.

### Research: CT as external support

All the pure mathematicians stated that their research interests are directed toward discovering relations among abstract theoretical concepts. Although most pure mathematicians talked about having used computations in one way or another as a tool to support this activity, they positioned the use of computers as external to the object of research:

> *“It helps you to calculate. It may be wrong to say that this is the primitive part of these machines, but in a way, it is. It is a tremendous help, but it is only a help.“*

Even though none of them positioned computers as part of their research goals, some still showed substantial engagement with computers. Some discussed spending considerable effort to work out the details of famous applications of theoretical results in widely used applications, such as algorithms for image compression, machine learning, and internet searches. In contrast, others talked about using computers and programming to verify examples and conjectures. However, such use was mostly backgrounded and either not seen as part of research at all, or as some type of “*extra bonus*”. One exception to this pattern was a mathematician who discussed using Macaulay. This software uses script-like programming to compute algebraic equations for certain geometrical objects beyond the reach of pen-and-paper computations, which he uses as part of his theoretical arguments. We do not place him among the numerical mathematicians, as he maintained a view of the purpose of his research to be theorems that make perfect sense in the absence of computers, and that his use of Macaulay left him ‘haunted’ by a feeling that part of the story was missing.

Pure mathematicians drew on real-world relevance mainly to motivate the significance of mathematics as a field, and did not see it as a constitutive part of their research. They still tended to be enthusiastic when talking about how ideas sometimes connect across disciplines, with one talking enthusiastically about how the notion of computers inspired theoretical work on modelling human cognition using neural networks already in the 1970s.

### Teaching: CT as marginalised

For the pure mathematicians, the use of CT mainly carried over from research when it was treated as motivational, provided conceptual support, or was overtly epistemically valuable. As an example, a mathematician who foregrounded his use of light-touch epistemic support of CT in research talked about using the same approach in his teaching, and that this allowed him to flip what he called the traditional "difficult first, easy later" sequence (theorem and proof first, examples and exercises later). He said this was well received by the students, and that they found it fun, perhaps simply because they, for once, were able to understand what they were supposed to do.

> *"I just wanted to create an assignment that would do something meaningful – where the programming part would lead to something…"*

Use of CT that pragmatically supported theoretical work or which was backgrounded in research tended not to be reproduced, unless there was a strong vocational imperative. In particular, there was little talk on the use of symbolic software in their teaching. Instead, these mathematicians emphasised the importance of developing students' analytical skills in their teaching. In the cases when use of computations was leveraged because of a feeling of institutional responsibility, such use was fragile and not seen as constitutive to the object of learning:

> *"And I was wondering… should I do some programming exercises […] And I was a bit lazy and didn't do it. I feel a bit guilty that I perhaps should have made more of it. On the other hand, the book isn't computationally oriented. Most of the exercises in the book can be done perfectly well without a program."*

At the level of thesis work, almost all mathematicians talked about allowing students the option to choose whether they wanted to use computations in their theses, regardless of their own engagement with computations. One said that he allowed students to draw on computations in thesis work, so that weaker students could get by without necessarily doing mathematical work.

## Numerical mathematicians: CT as internal but contested

We consider three of the interviewed mathematicians to be oriented toward research on numerical methods. We interpret their object of research as centred on the theoretical study of the computational properties of mathematical algorithms. These mathematicians tended to have a strained relationship with computations in their research, which carried over to their teaching.

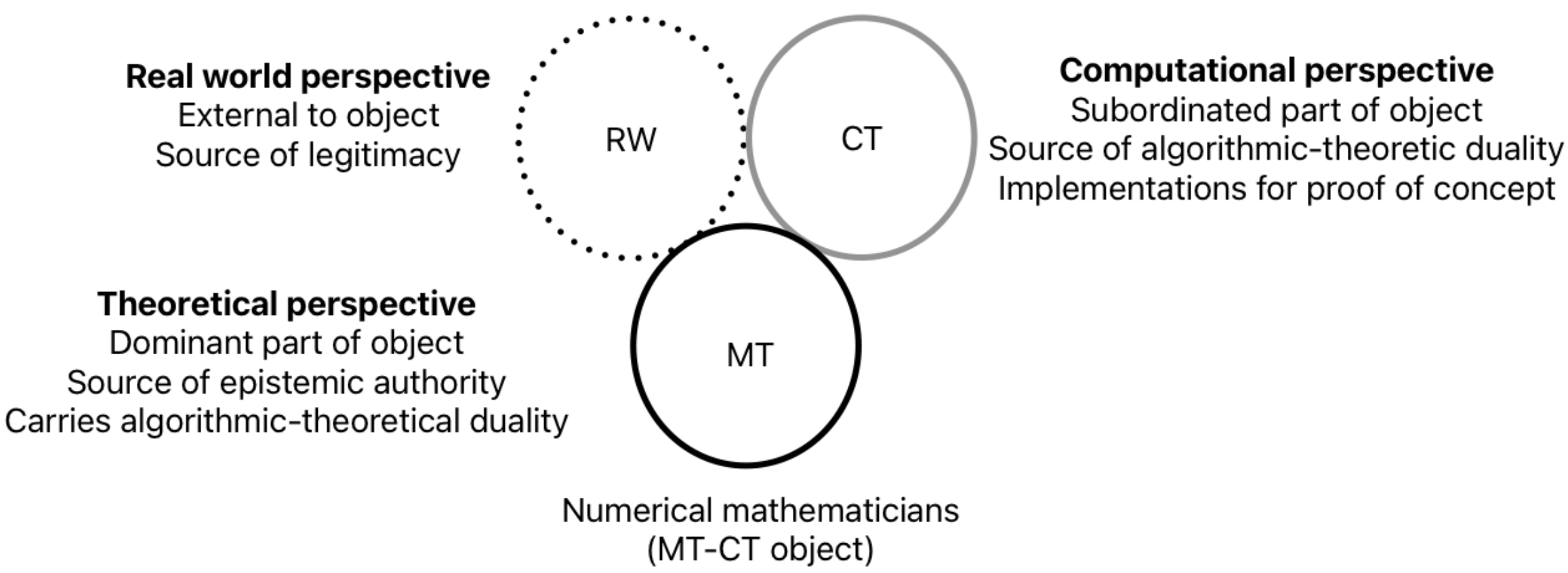

Figure 5: The numerical mathematicians conceptualised CT as part of their object, but it is subordinate to MT and considered primarily as a tool to illustrate and validate theoretical results, which were also considered to have a computational nature.

## Research: A computational-theoretical object of research

The numerical mathematicians explicitly described their object of research as having both theoretical and computational aspects. While they appeared to be skilled programmers who reported using programming in support of theoretical work, they came across as having a strained relation to programming. They most played down the role of computations in their mathematical work:

> *"The research work consists of developing numerical methods and then analysing those methods using functional analysis and other tools. And the actual programming of these methods – that's not the most important part, you could say."*

When asked whether they could dispense with the programming, they said that this would look bad as long as they were working on results that simulations could illustrate. However, as the following example shows, they struggled to explain why this was the case:

> *"It's kind of… It's a bit difficult to say. I don't know whether it's included simply because that's what people do – that it's standard practice in the field – or if... I would be a bit suspicious if an article about a numerical method didn't include graphs showing how it works in practice. I would find that a bit strange. I would take it as a sign that they might be trying to hide something – that it doesn't actually work in practice."*

We see this as evidence that their object of research differs from that of pure mathematicians. Indeed, should computers cease to exist, then the object of the pure mathematicians would remain intact (although perhaps harder to reach), the object of the numerical mathematicians would be fundamentally weakened. When one of them was asked whether the use of computers impacted his theoretical perspective on his work, he suspected that this may be so, but was not able to formulate how:

> *"So if you use a computer a lot, then you are choosing away something else. So it is at least the case that you may lose a bit of another perspective, I think. But it is probably not so easy to know what you have missed out on, in a way."*

Although they positioned programming as not part of their mathematical work, numerical mathematicians still saw their way of programming as distinct from that of computer scientists. In particular, one said that taking courses in 'programming' is not necessarily helpful if your goal is to learn 'algorithmic programming':

> *"But that is not very important for algorithmic programming. [...] it is somewhat different from the theory that they often learn in computer science courses."*

Notably absent from their interviews were mentions of collaboration with external partners and the use of real-world data, even though their terminology often draws on other disciplines, such as physics. When asked whether such use of terminology meant that their work had real-world applications, one replied: "*In principle there should be [a connection], but the road is long, I think*". Another talked about solving theoretical problems related to real-world problem applications handed to him by a collaborator who, he said, "*more or less just gave me the purely mathematical problems, which I solved for him."*

## Teaching: CT marginalised

The strained relationship of numerical mathematicians with programming in research also carried over to their teaching. For instance, a numerical mathematician who wanted to avoid programming in his research, talked about doing the same in teaching, even though he felt conflicted about it:

> *"You can choose to teach it purely theoretically – that is, a purely analytical treatment of the algorithms – and then that's it. Or you can teach it purely practically… just recipes and programming. Or something in between. I am probably more naturally inclined toward a purely analytical treatment of numerical methods. But I think students would benefit a lot – or do benefit a lot – from being forced to actually program the algorithms and make them work in practice."*

This pattern of tension across practices appeared in all three numerical mathematicians. What at first appears to be a mixed attitude toward programming in teaching is consistent with their research object. Most dismissed light programming use in 'pure' courses such as calculus or linear algebra, but all viewed it more positively when teaching aimed at insight into mathematical algorithms. Even the participant who used programming least recalled "aha" moments from programming as a student.

> *"When you have to program something, you really have to understand it, and you discover what you have misunderstood. And another thing is a bit more like… a bit harder to put your finger on, but I remember several times during my studies… when I was sitting there programming things and having those "aha" moments and rushes of joy from seeing that this actually works. This... in some magical way, I get something that resembles the exact solution."*

Despite recognising that computations provided some insight into algorithms, most preferred to avoid them in their teaching, stating that they found algorithmic programming pedagogically intractable in lower-level and graduate courses, and when supervising thesis students, even at the PhD level. However, one numerical mathematician was an exception to this. He said he found it fascinating to figure out how to have the computational and theoretical perspectives play off each other in his teaching. His use of code was most visible in his written course material and on assignments, and that most of his lectures are purely theoretical. His experience was that it could be difficult for students to coordinate the two perspectives.

## Applied mathematicians: CT as internal and stable

We consider six of the interviewed mathematicians to be applied mathematicians. We interpret their research objective as developing, implementing, and deploying algorithms to solve real-world problems. They all emphasised that their work had to be relevant to the real world, and that combining theory with computations was necessary to do this. However, they positioned CT in three different ways, each corresponding to a different positioning of CT in their teaching.

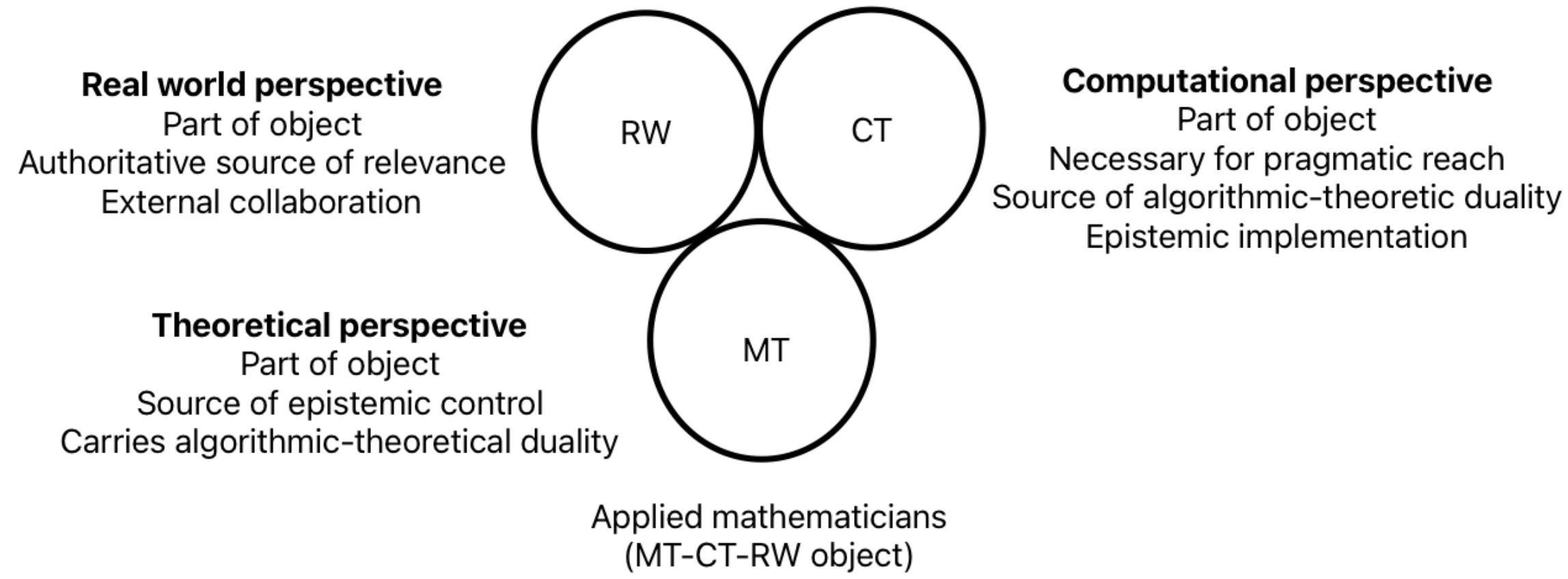

Figure 6: For the applied mathematicians, all three perspectives are positioned as part of their research objects. However, all seem to lean on some perspectives more than others.

### Research: Three different emphases of CT in object

The applied mathematicians all emphasised how the use of theory, computations, and collaboration with external partners is essential to their work. However, we find that they emphasised these perspectives in three different ways, which we now discuss.

Four of the interviewed mathematicians are real-world applied mathematicians with method emphasis (MT-CT). Their interest lay in developing efficient methods for solving mathematical models, similar to those of the numerical mathematicians. Interestingly, for them, efficiency is framed as a property of the algorithm and not of the implementation: *"The focus is really on the underlying methods. Others can implement them as efficiently as possible."* They emphasised real-world relevance and how they saw this as necessitating the use of computers. In the words of one of them, who started in pure mathematics: "*For a question to be important, to be useful, you need to be able to work with computers."* All four applied mathematicians with a method emphasis reported high engagement with AI and machine learning as research objects. They seemed to view AI/ML as a RW-relevant but black-boxed model, which they saw as their job to study.

> *But that we actually have a role in serving society, and that in a world now with artificial intelligence on all sides, mathematics departments and mathematicians just have to* *be extremely grateful for that. There won't be any unemployed mathematicians anymore. And I believe they also need to be able to code.*

One of the mathematicians placed strong emphasis on modelling (MT-RW). His research interests lie in the relation between mathematical models and real-world phenomena, with the role of implementation being backgrounded. This mathematician argued that, in the end, all results must be checked against real-world data and experiments, since both the modelling and the math could be wrong. While he acknowledged that computers are "part of the game" and necessary to generate model outputs he can compare to real-world data, he said he had no interest in spending more time than necessary on implementation. In particular, he prefers to code from scratch, avoiding any type of black-boxing.

> *"Programs are increasingly functioning as black boxes. And it's by no means certain that these so-called wave researchers have written their own programs either. They sit there using them as black boxes that someone else. That almost makes me worried, because I feel that you miss out on some of the deeper analytical insight that I've been used to having."*

The final applied mathematician had a strong emphasis on developing and deploying computational solutions to real-world problems (CT-RW). He was the only mathematician interviewed who explicitly placed equal emphasis on mathematics and computations as tools: "*I wouldn't say that one tool is much better than another. What matters is solving the damn problem. Whether mathematics is harder than programming or physics or whatever – I don't know."* Specifically, he saw mathematics as being about modelling and debugging, while computers provide pragmatic reach. However, he was also the mathematician who emphasised implementation the most. When asked how he handles bugs in complicated projects, he pointed to the use of coding assistants and how they are fundamentally changing the way people program:

> *Now, when I program [...] it's basically a matter of Googling something, or I might even use Copilot or something like that to set up a rough sketch, and then I go from there. So programming has perhaps changed.*

## Teaching: Real-world relevance as a significant mediator

The three ways applied mathematicians position CT correspond to different ways computation carries over into teaching.

The method-oriented mathematicians, with MT-CT emphasis, tended to talk about CT as part of their pedagogy, even in foundational courses. They described using computations to help students explore theoretical concepts through simulations to gain hands-on experience and familiarise themselves with issues that they say only become visible during implementation. For them, the RW perspective is important. They talked about using computations pedagogically to allow for earlier exposure to key ideas since students who only take their first course could also end up in a job as a 'number cruncher'. We see this sense of societal responsibility as instrumental to their decision to engage with machine learning and generative AI in their teaching. They said that it is exactly because these methods were controversial that students needed to encounter them during their studies.

The model- and implementation-oriented mathematicians marginalised CT in foundational courses, but for different reasons. The model oriented mathematician, with MT-RW emphasis, said that while he supported that students learn computations early, he did not see it as defensible for him to be part of it: "*I do not think we can afford to add more programming there. Because then I believe it would come at the expense of the mathematics they ought to learn*." He insisted on computations not be used by students as black-boxes. The implementation-oriented mathematician, with CT-RW emphasis, did not express any concerns over black-boxed use of computations. For him, there was no point in complicating courses such as calculus with computations, as students would just get bogged down debugging code for simple problems they could just solve with pen-and-paper: *"It is no longer elegant or instructive; it just becomes tedious."* However, he admitted he had not given this much thought, suggesting that his pragmatic approach to mathematics may have been the underlying reason for his low engagement.

On the level of thesis work, all applied mathematicians talked about giving students the freedom to decide whether their projects should lean toward theory, modelling or implementation. They emphasised giving students the freedom to follow their interests. One said, "*So I give them full freedom to work within precisely the area or areas they themselves prefer. I do not think I have ever required a student to write code against their will*."

# Discussion

In this study, we used CHAT and CoP to analyse interviews with 15 mathematicians about how they position computations across their research and teaching activities.

Across the sample, how mathematicians configured the research object carried over to teaching. Pure mathematicians treated CT as external support and introduced it mainly as motivation, epistemic scaffolding, or vocational accommodation. Numerical mathematicians treated CT as internal to work on algorithms, yet still subordinated it to theory, producing an ambivalent stance in teaching. Applied mathematicians treated CT as necessary to mathematically relevant work in the world, but differed in whether they foregrounded method, model, or implementation. Across all groups, foundational teaching foregrounded MT as epistemic control more strongly than advanced teaching and thesis supervision. These patterned differences suggest that perceived tensions between CT and MT depend not only on methods, but on how the object of mathematical work is framed.

We now discuss how these positionings help answer the research questions posed in the Background section.

**RQ1:** *How do mathematicians position computational thinking relative to mathematical thinking within their research activity?*

Real-world relevance strongly shapes how mathematicians position CT relative to MT. We interpret this as an ontological shift: When RW was external, theoretical relations dominated, and CT remained supportive. When RW became constitutive, theoretical concepts were more readily treated as algorithmic and computational, making CT part of the research object rather than merely a tool.

**RQ2:** *How does this positioning mediate engagement with computational thinking in teaching?*

Mathematicians' ontological positioning remains largely stable across research and teaching, but varies by educational level. In foundational teaching, analytic theory is more authoritative, and CT is often suppressed; in advanced teaching and thesis work, computation is more often accommodated for epistemic, research, or vocational reasons.

**RQ3:** *How does this positioning mediate mathematicians' talk about tensions or synergies between mathematical and computational thinking in education?*

Mathematicians' talk about tensions and synergies depends strongly on their ontological positioning. Although some described epistemic uses of computation as an on-ramp to theoretical insight, those who strongly emphasised theory as a source of analytic control tended to frame computation as a distraction and analytic competence as a prerequisite for its use. Where RW was constitutive, CT was more often framed as synergistic, though mainly for algorithmic rather than purely analytic understanding. The resulting tensions concerned not only epistemic incompatibility, but also the pedagogical cost of coordinating programming and mathematics.

## Theoretical contributions and implications for education

Our study suggests that tensions and synergies between mathematical and computational thinking are not only matters of epistemological compatibility, but also reflect variation in how mathematicians ontologically position computation within mathematical work. Tensions occur when this ontological positioning comes into conflict with actual practices (their own or others). The findings of Huang et al. (2021) seem to add a further layer, as movement between contexts may shape such positionings.

The need to attend to variation in tacit ontological considerations limits consensus-oriented approaches to computational thinking in mathematics, such as that of Kallia et al. (2021), since an emphasis on common denominators may obscure meaningful distinctions in how programming is positioned. In particular, numerical and applied mathematicians tend to regard their use of programming as distinct from the way programming is conceptualised in computer science. This suggests a tension with Wing's (2006) framing of CT as a set of transferable skills into other domains, and that taxonomies for use of programming grounded in computer science (e.g., Broley et al., 2018; Weintrop et al., 2016) may obscure distinctions that are meaningful within mathematics.

The three-perspectives model proposed in this study extends the epistemic matrix from Olsen et al. (2016). Whereas the epistemic matrix classified mathematicians according to both object *and* tool positionings, the three-perspectives model frames mathematics as a theory-computation-model dialectic. That is, mathematicians gain complementary insight into mathematical concepts by studying them in terms of theoretical relations, by their computational nature, and by how they model real-world phenomena. This dialectic is a cause of tension as mathematicians do not agree on where to draw a line between what counts as mathematics and what does not.

While our observations on mathematicians' talk about machine learning and AI should be taken as tentative, they show a strong correlation with how the mathematicians position themselves with respect to the three-perspectives model. This suggests that the model could serve as an analytic tool for making sense of mathematicians' engagement with emerging technologies, even before stable practices have formed.

In education, teachers and students have to negotiate different disciplinary contexts, which may shift ontological views (Huang et al., 2021). Our findings suggest structured ways to navigate these contexts to find meaningful and therefore stable ways to integrate computational and mathematical thinking. In particular, in a context that emphasises theoretical insight, use of programming should be backgrounded and be clearly in support of theoretical insight. To justify the added effort of coordinating mathematical and computational thinking, computation should be tied to curricular goals external to mathematics that highlight its pragmatic reach. In such settings, students can use computation to mobilise mathematical thinking in order to analyse data, explore models, or understand phenomena. Here, CT may become meaningful not because it is simply added to mathematics, but because it extends what students can do with mathematics.

# Limitations

In accordance with our interpretive and socio-cultural theoretical framing, we acknowledge that our results depend on context and the co-construction of meaning in encounters between the team of authors and the interviewed mathematicians. In particular, we deliberately chose a critical site because it had the special property that integration of computations in teaching has strong institutional support. However, we did not observe the practices of the mathematicians, and instead considered their discourse on their practices as our unit of analysis, using information available online to triangulate our interpretations.

While allowing us to address methodological concerns by Weber et al. (2020), the use of a researcher with experience in mathematical research to conduct interviews and analyses introduces risks associated with an insider perspective. This risk has been mitigated by this research being part of a long-ongoing research effort, where analyses of the same data set have been discussed not only within the team of authors, but also with other experts in mathematical and education research, but informally, at conferences and seminars.

The interviews were conducted in the fall of 2023, which means that the development of generative AI was still at a relatively early stage. This means both that the nature of the mathematicians' views on the interface between computational thinking and mathematics may have changed since then, and, in particular, how this relates to AI and machine learning. However, rather than a limitation, we see this as providing interesting insight into how the mathematicians talked about programming at a critical point in time, and which will provide an interesting baseline for how such views - and the stability of positioning of CT relative to their research objects - evolve since then.

## Future research

Our study suggests several directions for further research. First, it is necessary to conduct further empirical studies on mathematicians to validate our proposed three-perspectives model. Second, it would be interesting to explore whether extensions of the model can be formulated to explain the positioning of CT and MT relative to disciplinary methodologies of other domains, such as physics or teacher education. Third, it would be interesting to explore whether the model can be extended in the direction of students' perspectives, to better understand how they position MT and CT relative to *their* object of learning. Fourth, while we suggest some implications for education, these should be interpreted, and any educational implementations warrant independent research. Fifth, we observe a tension between how CT is conceptualised by computer scientists and mathematicians. It would be interesting to revisit taxonomies related to the use of CT in the mathematics concept while taking ontological positionings into account. Finally, we observe that patterns of engagement with generative AI and machine learning across teaching and research practices seem to resonate with how CT and MT are positioned as objects of research. Further research is necessary to confirm that this pattern continues to hold and what the implications are for education.

# Declarations

The authors declare that they have no competing interests.